\documentclass[10pt]{article}
\usepackage[english]{babel}
\usepackage{amsfonts}
\usepackage[dvips]{graphicx}

\begin{document}

\title{\bf Searching the (really) real general solution of 2D Laplace differential equation}

\date{2009, October}

\author{\bf Gianluca Argentini \\
\normalsize{[0,1]Bending - Italy}\\
\normalsize 01bending@gmail.com \\
\normalsize gianluca.argentini@gmail.com \\}

\maketitle

\noindent{\bf Abstract}\\
This is not a new result. Purpose of this work is to describe a method to search the analytical expression of the general real solution of the two-dimensional Laplace differential equation. This thing is not easy to find in scientific literature and, if present, often it is justified with the assertion that an arbitrary analytic complex function is a solution of Laplace equation, so introducing the condition of complex-differentiability which is not really necessary for the existence of a real solution. The question of the knowledge of real exact solutions to Laplace equation is of great importance in science and engineering.\\


Consider the 2D Laplace equation $\partial_{xx} u + \partial_{yy}u = 0$ for a function $u=u(x,y)$ defined on a region $\Omega \subset \mathbb{R}^2$. If one, e.g. for physical or engineering reasons, is interested in the symbolic (analytical) expression of the general {\it real} solution, usual technical good literature (see \cite{polyanin} for a general reference) gives three categories of answers:\\

\noindent {\it a.} theorems of existence and uniqueness for boundary value problems, very important results but often without a symbolic expression of the solution;\\

\noindent {\it b.} applications of the method of separable variables, which often gives a symbolic solution of the differential equation but only in the restricted range of functions of the form $F(x)G(y)$;\\

\noindent {\it c.} an information, unfortunately often not completely justified, about the fact that an {\it analytic} (in the sense holomorphic) complex function $F(z)$ is harmonic, that is $\Delta F(z)=0$, where $0$ is the origin of the complex plane.\\

\noindent In particular, in the latter case the function $F$ must be $\mathbb{C}$-derivable, that is $\partial_x F = -i\partial_y F$, which is a condition very stronger than the differentiability of the real and imaginary part; therefore, complex analytic functions are a restricted range of all the possible solutions of the Laplace equation.

\noindent There is another way to search the form of the solution: the use of a symbolic mathematical software. If we ask {\it Mathematica} to solve $\partial_{xx} u + \partial_{yy}u = 0$, the built-in function {\ttfamily DSolve} answers

\begin{equation}\label{genSolutionNotShown}
	u(x,y) = F(y + jx) + G(y - jx)
\end{equation}

\noindent where $j^2=-1$ and $F$, $G$ are arbitrary functions. The relative tutorial (\cite{DSolve}) gives no more informations about the nature of these {\it arbitrary} functions. Note that, in general, $F(y+jx)$ and $G(y-jx)$ are complex numbers. Previous formula (\ref{genSolutionNotShown}) is present in (\cite{sneddon}), but only as exercise to show that it is a solution of Laplace equation.

\begin{figure}[h!]\label{DSolve}
	\begin{center}
	\includegraphics[width=12cm]{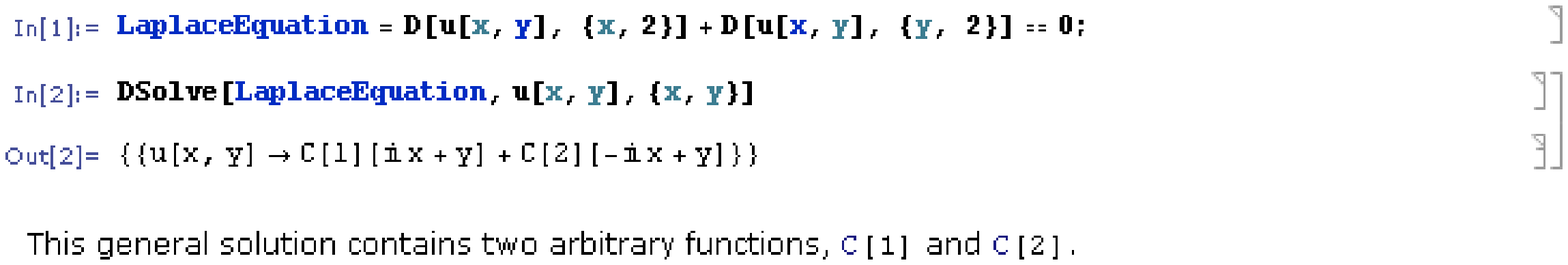}
	\caption{{\it Output from} {\textnormal {\ttfamily DSolve}}. }
	\end{center}
\end{figure}

\noindent Hence, what is the mathematical form of the real general solution of the Laplace equation? In this work I try to find it with a method which is not standard in the usual literature.\\

\noindent Consider the laplacian operator $\Delta = \partial_{xx} + \partial_{yy}$ as an algebraic object defined on the set of twice partially differentiable real functions defined on $\mathbb{R}^2$. We can define the object $(\partial_x)^2 = \partial_{xx}$, that is $(\partial_x)^2 f = \partial_x \partial_x f$ for every $f \in C^2(\mathbb{R}^2,\mathbb{R})$. Then it can be easily stated the following identity:

\begin{equation}\label{identityLaplacian}
	\Delta = (\partial_x + j\partial_y)(\partial_x - j\partial_y)
\end{equation}

\noindent because on $C^2(\mathbb{R}^2,\mathbb{R})$ the two single elementary operators $\partial_x$ and $\partial_y$ commute: $(\partial_x)(j \partial_y) f = (j \partial_y)(\partial_x) f$. Define $D_j = (\partial_x + j\partial_y)$ and $D_{-j} = (\partial_x - j\partial_y)$. Two similar operators are defined, in the discussion of the Cauchy-Riemann equations, in the first edition of (\cite{rudin}).\\

\noindent Let be $u \in C^2(\mathbb{R}^2,\mathbb{R})$; then, being $D_{-j}D_j u = D_jD_{-j} u$, $u$ is a solution of $\Delta u = 0$ if $D_j u = 0$ or $D_{-j} u = 0$. Consider the first case:

\begin{equation}\label{Dj}
	\partial_x u + j\partial_y u = 0
\end{equation}

\noindent It is easy to apply the method of characteristics (see \cite{jeffrey}) for its symbolic resolution. The ordinary differential equation resolved by the characteristics is

\begin{equation}
	y'(x) = j
\end{equation}

\noindent which is satisfied by the family $y = jx + c$, with $c$ arbitrary (complex) constant. Along each characteristics the equation (\ref{Dj}) becomes (note that the complex function $jx+c$ is $\mathbb{R}$-derivable with respect to the real variable $x$)

\begin{equation}
	d_x \left( u(x, jx + c) \right) = 0
\end{equation}

\noindent so for every choice of $c$ there is an arbitrary (complex) constant $G$, depending on $c$, such that $u(x, jx + c) = G(c)$. But for every value of $y$ there is a characteristics, namely the line corresponding to the constant $c = y - jx$, passing for the point $(0,y)$ of $\mathbb{R}^2$), so if we take a (complex) function $G$ whose real and imaginary part are twice $\mathbb{R}$-differentiable on $\mathbb{R}^2$, the general real solution of (\ref{Dj}) is

\begin{equation}\label{solDj}
	u_1 = \mathbb{R}e \left[ G(y-jx) \right]
\end{equation}

\noindent In the same manner, one finds that $u_2 = F(y+jx)$ is the general complex solution of the second case, $D_{-j} u = 0$.\\
\noindent Now, suppose $\Delta u = D_jD_{-j} u = 0$. If $D_{-j} u = 0$ or $D_{j} u = 0$, the case has already been treated. Then, suppose e.g. $D_{-j} u \neq 0$. Therefore $D_j \left( D_{-j}u \right) = 0$, so from (\ref{solDj}) there is a function $G$, with real and imaginary part twice $\mathbb{R}$-differentiable on $\mathbb{R}^2$, such that $D_{-j}u = G(y-jx)$. Now, to find $u$, we have to resolve the non homogeneous PDE

\begin{equation}\label{pdeNonHom}
	\partial_x u - j\partial_y u = G(y-jx)
\end{equation}

\noindent But the homogeneous equation $D_{-j} u = 0$ has general solution $u_h = F(y+jx)$. We try to find a particular solution $u_p$ of (\ref{pdeNonHom}). Let be $u_p(x,y) = a(x,y)+jb(x,y)$. Then

\begin{equation}
	D_{-j} u_p = (\partial_x a + \partial_y b) + j (-\partial_y a + \partial_x b)
\end{equation}

\noindent and the two following first order pdes must hold ($G_{\mathbb{R}} = \mathbb{R}e[G]$, $G_{\mathbb{I}} = \mathbb{I}m [G]$):

\begin{equation}
	\partial_x a + \partial_y b = G_{\mathbb{R}}, \hspace{0.1cm} -\partial_y a + \partial_x b = G_{\mathbb{I}}
\end{equation}

\noindent A possible way to find a solution is the choice

\begin{eqnarray}
	\partial_x a = \frac{1}{2}G_{\mathbb{R}}, \hspace{0.1cm} -\partial_y a = \frac{1}{2}G_{\mathbb{I}} \nonumber \\
	\partial_x b = \frac{1}{2}G_{\mathbb{I}}, \hspace{0.1cm} \partial_y b = \frac{1}{2}G_{\mathbb{R}} \nonumber
\end{eqnarray}

\noindent which can be resolved by quadrature. Let be $\mathbb{A}$ the complex function so defined:

\begin{equation}
	\mathbb{A}(z) = \mathbb{A}(x+jy) = a(-y,x)+jb(-y,x)
\end{equation}

\noindent Note that $\mathbb{A}(y-jx) = a(x,y)+jb(x,y)$, so the function $u_p(x,y)=\mathbb{A}(y-jx)$ is a particular solution of (\ref{pdeNonHom}). We have shown that if $\Delta u = 0$, then $u$ has the form (\ref{genSolutionNotShown}).\\

\noindent In conclusion, the general real solution of the Laplace equation is

\begin{equation}\label{generalRealSolutionLaplace}
	u(x,y) = \mathbb{R}e \left[ F(y+jx) + G(y-jx) \right]
\end{equation}

\noindent Note that for our considerations it is not necessary the $\mathbb{C}$-differentiability of $F$ and $G$. Therefore, the assertion about Laplace equation, made by {\it Mathematica} and by some texts on complex analysis, is now more clear when transposed in the context of real functions:\\

\noindent {\it the general real solution of Laplace equation $\partial_{xx}u + \partial_{yy}u = 0$ is $u(x,y) = \mathbb{R}e \left[ F(y+jx) + G(y-jx) \right]$ where $F$ and $G$ are arbitrary complex functions such that the real part of $F(y+jx) + G(y-jx)$, considered as function of $(x,y)$, is twice $\mathbb{R}$-differentiable}.\\

\noindent {\it Example}. Suppose we want to find a solution of the differential problem $\partial_{xx}\Psi + \partial_{yy}\Psi = 0$ on the domain $\Omega=\{ (x,y): x \geq 0, 0 \leq y \leq h \}$, where $h > 0$, with boundary condition $\Psi = 0$ on $\partial \Omega$.\\
\noindent Such a problem arises in fluid dynamics when $\Psi$ is the stream function of a flow in $\Omega$, semi-infinite bay, with velocity field $(\partial_y \Psi, -\partial_x \Psi)$ (see \cite{madani}).\\
\noindent Note that the problem has not a unique solution. The trivial function $\Psi = 0$ is solution, and for every $\Psi$ solution of the problem, $k\Psi$ is solution too for every constant $k$. We try to construct a not trivial solution from the general formula (\ref{generalRealSolutionLaplace}). Being $\Psi(0,y)=0$ for every $y \in \left[0,h\right]$, a possible solution could be $\Psi = \sin(x)f(y)$ with $f$ such that $f(0)=f(h)=0$. One knows that for the complex function $\cos(jx+y)$ the following identity holds:

\begin{equation}
	\cos(jx+y) = \cosh(x)\cos(y)-j\sinh(x)\sin(y)
\end{equation}

\noindent The function $j\cos(jx+y)$ is of type $F(y+jx)$, and its real part is a twice $\mathbb{R}$-differentiable function. For the boundary conditions, note that it is sufficient to apply a simple change of variable $jx+y \rightarrow \frac{n\pi}{h}(jx+y)$, where $n \in \mathbb{N}$. The real part becomes

\begin{equation}
	\Psi(x,y) = \sinh\left( \frac{n\pi x}{h} \right) \sin\left( \frac{n\pi y}{h} \right)
\end{equation}

\noindent that is a real solution of the boundary problem. The same problem, with the same solution, is treated in (\cite{madani}), where separation of variables is applied.\\

\noindent {\it Application}. Find the general real solution of the PDE

\begin{equation}\label{stressEquation}
	\partial_{xx}U - \partial_{yy}U = 0
\end{equation}

\noindent This equation appears in the description of the shear stress of a plane flow with velocity field $(\partial_y U, -\partial_x U)$. We can repeat previous argument with $\partial_{xx} - \partial_{yy} = (\partial_x + \partial_y)(\partial_x - \partial_y)$. But it is more interesting if we apply the variables transformation $X = x$, $Y = jy$, where $j^2=-1$. Then (\ref{stressEquation}) becomes $\partial_{XX}U + \partial_{YY}U = 0$ where $U=U(x(X),y(Y))$, which is the previous Laplace equation. Its general real solution is $U(x,y) = \mathbb{R}e \left[ F(Y(y)+jX(x)) + G(Y(y)-jX(x)) \right]$, therefore, being $Y+jX=j(y+x)$ and $Y-jX=j(y-x)$, we can write

\begin{equation}
	U(x,y) = \mathbb{R}e \left[ F(y+x) + G(y-x) \right]
\end{equation}

\begin{figure}[ht!]
	\includegraphics[width=2cm]{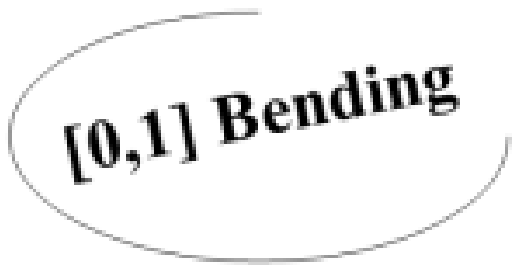}
\end{figure}
\noindent \tiny {\bf [0,1]Bending}, a Design Studio in Italy dedicated to computational engineering for scientific and industrial applications.\\
\noindent gianluca.argentini@gmail.it


\begin{thebibliography}{9}

\bibitem{DSolve} {\it Differential Equation Solving with {\textnormal {\ttfamily DSolve}}}, Mathematica 7.1, www.wolfram.com, 2009

\bibitem{jeffrey} A.Jeffrey, {\it Applied Partial Differential Equations: An Introduction}, Academic Press, 2002

\bibitem{madani} R.Malek-Madani, {\it Advanced Engineering Mathematics}, vol. $\mathbf{2}$, Addison-Wesley, 1998

\bibitem{polyanin} A.Polyanin, {\it Handbook of linear partial differential equations for engineers and scientists}, Chapman \& All/CRC, 2002

\bibitem{rudin} W.Rudin, {\it Real and Complex Analysis}, McGraw-Hill, 1966

\bibitem{sneddon} I.Sneddon, {\it Elements of Partial Differential Equations}, Dover, 2006

\end{thebibliography}
\end{document}